\documentclass[
submission
]{dmtcs-episciences}

\usepackage[utf8]{inputenc}
\usepackage{subfigure}
\usepackage{mathrsfs}
\usepackage{amssymb}
\usepackage{mathrsfs}
\usepackage{enumerate}
\usepackage{hhline}
\usepackage{amsmath}
\usepackage{color}
\usepackage{graphicx}
 \usepackage{epstopdf}
 \usepackage{lineno,hyperref}
 \usepackage{natbib} 
 \setcitestyle{numbers,square}
  
\newtheorem{proposition}{Proposition}[section]
\newtheorem{pro}[proposition]{Problem}%[section]
\newtheorem{lemma}[proposition]{Lemma}%[section]
\newtheorem{theorem}[proposition]{Theorem}%[section]
\def\q{\hspace*{\fill}$\Box$\medskip}

\author{Sarula Chang\affiliationmark{1}\thanks{Supported by Program for improving the Scientific Research Ability of Youth Teachers of Inner Mongolia Agricultural University(No.BR230110); NSF of China (Nos.12361069, \,12171089,\,12271235); NSF of
Fujian (No.2021J02048); Research Fund of Xiamen university of technology(Nos.YKJ20018R, XPDKT20039.)}
  \and Jianxi Li\affiliationmark{2}
  \and Yirong Zheng\affiliationmark{3}}
\title[A characterization on trees $T$ with $m(T, \lambda)=p(T)-2$]{A characterization on trees $T$ with $m(T, \lambda)=p(T)-2$}
\affiliation{
  College of Science, Inner Mongolia Agricultural University, Hohhot, Inner Mongolia,  China\\
  School of Mathematics and Statistics, Minnan Normal University, Zhangzhou, Fujian,  China\\
 School of Mathematics and Statistics, Xiamen University of Technology, Xiamen, Fujian,  China}
\received{2000-00-00}
\revised{2000-00-00, 2000-00-00, 2000-00-00}
\accepted{2000-00-00}
\begin{document}
\publicationdetails{VOL}{2015}{ISS}{NUM}{SUBM}
\maketitle
\begin{abstract}
  Let $m(G,\lambda)$ be the multiplicity of an eigenvalue $\lambda$ of
a connected graph $G$. Wang et al. [Linear Algebra Appl. 584(2020),
257-266] proved that for any connected graph $G\neq C_n$,  $m(G,
\lambda) \leq 2c(G) + p(G) -1$, where $c (G) = |E(G)| - |V (G)| + 1$
and $p(G)$  are the cyclomatic number and the number of pendant
vertices of $G$, respectively.  In the same paper, they proposed the
problem to characterize all connected graphs $G$ with
 eigenvalue $\lambda$ such that $m(G, \lambda) =2c (G)+ p(G)-1$.
Wong et al. [Discrete Math. 347(2024), 113845] solved this problem
for the case when $G$ is a tree by characterizing all trees $T$ with
eigenvalue $\lambda$ such that $m(T , \lambda) = p(T )-1$. In this
paper, we further provide the structural characterization on trees
$T$ with eigenvalue $\lambda$ such that $m(T , \lambda) = p(T )-2$.
\end{abstract}

\section{Introduction}
\label{sec:in}
All graphs considered in this paper are finite, undirected and
simple. Given a graph $G=(V(G), E(G))$ with vertex set $V(G) = \{v_1, \ldots, v_n\}$, the \emph{adjacency matrix} of $G$ is the $n\times n$ matrix $A(G) = (a_{ ij })$ such that $a _{ij}=1$
 if vertices $v _i$ and $v_ j$ are adjacent; and $a _{ij} =0$, otherwise.
 The eigenvalues of $A(G)$ are called the \emph{eigenvalues} of $G$. By $m(G, \lambda)$ we denote the \emph{multiplicity} of a real number $\lambda$ as an eigenvalue of $G$.
% As usual, we call $|E(G)| -|V (G)| + \omega(G)$ the \emph{cyclomatic number} (or \emph{dimension of cycle space})  of $G$ and  denote it by $c(G)$, where $\omega(G)$ is the number of connected components of $G$. In particular, for a connected graph $G$, if $c(G) = 0$, then $G$ is a \emph{tree}.
Let $c(G)=|E(G)| -|V (G)| + \omega(G)$ be the \emph{cyclomatic
number} (or \emph{dimension of cycle space}) of $G$, where
$\omega(G)$  is the number of connected components of $G$. In 2020,
Wang et al. \cite{Wang1} gave an upper bound for $m(G, \lambda)$ in
terms of $c(G)$ and $p(G)$, where  $p(G)$ is the number of pendant
vertices in $G$.

\begin{theorem}[\cite{Wang1}]\label{1.1}
Let $G$ be a connected graph with at least two vertices.
If $G$ is not a cycle, then $m(G, \lambda)\leq 2c (G) +p(G) -1$ for any $\lambda \in R$.
\end{theorem}

In the same paper, they proposed the following problem for connected
graphs $G$ ($\omega(G) =1$):

\begin{pro}\label{p1-2}
Characterize all connected graphs $G$ with eigenvalue $\lambda$ such
that $m(G, \lambda) = 2c (G) + p(G) - 1$.
%with $m(G, \lambda)\leq 2c (G) +p(G) -1$ for any $\lambda \in R$.
\end{pro}

{\bf Remark:} Indeed, Problem~\ref{p1-2} has been solved for two
special cases, the case when $\lambda = 0$ and the case when
$\lambda = -1$. Chang et al. \cite{Chang} characterized the
connected graphs $G$ without pendant vertices which has 0 as an
eigenvalue of multiplicity $2c (G) + p(G) -1$ and Wang et al.
\cite{Wang} and Chang et al. \cite{Chang1} gave a complete
characterization for $G$ with $m(G, 0) = 2c(G) + p(G) -1$,
respectively. Recently, Zhou et al. \cite{Zhou} characterized the
connected graphs $G$ with $-1$ as an eigenvalue of multiplicity $2c
(G) + p(G) - 1$. Before these works, the authors in \cite{Ma} have
characterized the connected graphs $G$  which has 0 as an eigenvalue
of multiplicity (i.e. nullity) $2c (G) + p(G)$.

In particular, when $G$ is a tree ($c(G)=0$), then Theorem~\ref{1.1}
becomes that for any tree $T$ with eigenvalue $\lambda$, $m(T ,
\lambda) \leq p(T ) - 1$. Recently, Wong et al.\cite{Wong}  give the
following complete characterization for trees $T$ with $m(T,
\lambda) = p(T) -1$,  where $\lambda$ is an eigenvalue of $T$, which
solves Problem~\ref{p1-2} when $G$ is a tree.

%A corollary follows immediately.
%\begin{corollary}[\cite{Wang1}]\label{1.2}
% Let $T$ be a tree and $\lambda$ be a real number, then $m(T , \lambda) \leq p(T ) - 1$.
%\end{corollary}

%In \cite{Wang1} the authors concluded with the following remark:

%The following problem seems difficult, however,
%interesting: To characterize the connected graphs $G$ together with
%the eigenvalue $\lambda$ such that $m(G, \lambda) = 2c (G) + p(G) - 1$.

\begin{theorem}[\cite{Wong}]\label{1.3}
Let $T$ be a tree with $\lambda$ as an eigenvalue of $T$. Then $m(T , \lambda) = p(T ) - 1$ if and only if
\begin{enumerate}[(1)]
\item there are two co-prime positive integers $i$ and $m + 1$ with $1 \leq i \leq m$ such that $\lambda= 2cos \frac{i\pi }{m+1}$;
\item $T\in\Gamma_k(\lambda)$, where $k =\gamma (T )$ is the number of major vertices and $\Gamma_k(\lambda)$ is defined as below.
    \end{enumerate}
\end{theorem}

In \cite{Wong}, $\Gamma_j(\lambda)$ is defined as below  with respect to a real number
 $\lambda = 2\cos\left(\frac{i\pi}{m+1}\right)$, where $i\in [1, m]$ and $m + 1$ are co-prime positive integers.
\begin{itemize}
\item Denote by  $\Gamma_{0}(\lambda) $ the set consisting of the paths  $P_n$ with $n \equiv m (mod ~m +1)$, or equivalently, $n +1$ is a multiple of $m +1$. Namely, $\Gamma_{0}(\lambda) = \{P_n|n \equiv m (mod~ m + 1)\}$.

\item  By  $\Gamma_{1}(\lambda) $, we denote the set of all trees $T$ satisfying the following conditions: Each $T\in\Gamma_1(\lambda)$ has a unique major vertex $w_1$  such that $T - w_1$ is the union of at least three components from $\Gamma_{0}(\lambda) $, where  the neighborhood of $w_1$ in each component is a pendant vertex of that component.

\item  By  $\Gamma_{2}(\lambda) $, we denote the set of all trees $T$ satisfying the following conditions: For each $T\in\Gamma_2(\lambda)$, $\gamma (T ) = 2$, and there is a major vertex $w_2$ of  $T$ such that $T - w_2$ has exactly one component from $\Gamma_{1}(\lambda) $ and all other components come from $\Gamma_{0}(\lambda) $, where  the neighborhood of $w_2$ in each component is a pendant vertex of that component.

 \end{itemize}

    Based on the definition for $\Gamma_{j-1}(\lambda) $ , we can define  $\Gamma_{j}(\lambda)$, recursively.

\begin{itemize}
\item  By   $\Gamma_{j}(\lambda) $, we denote the set of all trees $T$ satisfying the following conditions: For each $T\in\Gamma_j(\lambda)$, $\gamma (T ) = j$, and there is a major vertex $w_j$ of  $T$ such that $T - w_j$ has exactly one component from $\Gamma_{j-1}(\lambda) $ and all other components come from $\Gamma_{0}(\lambda) $, where  the neighborhood of $w_j$ in each component is a pendant vertex of that component.
 \end{itemize}

\begin{figure}[h]
\centering
\includegraphics[width=0.9\textwidth]{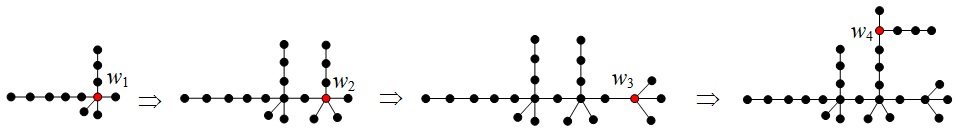}
\caption{ A flow chart on how to construct a tree $T$ in $\Gamma_{4}(0)$}\label{Fig.1}
\end{figure}

Fig.~\ref{Fig.1} is a flow chart explaining how to construct a tree
in $\Gamma_{4}(0) $, where $0 = 2\cos \frac{\pi}{2}$.

Along this line, in this paper, we further characterize trees $T$
with eigenvalue $\lambda$ such that $m(T , \lambda) = p(T ) - 2$.
Our results can be read as follows.

\begin{theorem}\label{1.40}
Let $T$ be a tree with $p(T ) =3 $ and $\lambda$ be an eigenvalue of
$T$. Then $m(T , \lambda) = p(T ) - 2=1$ if and only if
$T\in\Gamma^{2}_{1}(\lambda)$, where $\Gamma^{2}_{1}(\lambda)$ is
defined as below.
\end{theorem}

\begin{theorem}\label{1.4}
Let $T$ be a tree with $p(T ) \geq4 $ and $\lambda$ be an eigenvalue
of  $T$. Then $m(T , \lambda) = p(T ) - 2$ if and only if
$T\in\Gamma^{2}_{k}(\lambda)$, where $k =\gamma (T )$  is the number
of major vertices, where  $\Gamma^{2}_{k}(\lambda)$  is defined as
below.
\end{theorem}

We   construct a class of trees $\Gamma^{2}_{j}(\lambda)$ with
respect to a real number $\lambda = 2\cos \frac{i\pi}{m+1}$, where
$i\in [1, m]$ and $m + 1$ are co-prime positive integers.

 \begin{itemize}
 \item Denote by  $\Gamma^{2}_{0}(\lambda)$ the set  consisting of trees $T\in \Gamma_{0}(\lambda)$ by deleting a  pendant
 vertex.

\item  By  $\Gamma^{2}_{1}(\lambda)$, we denote the set of all
trees $T$ satisfying the following conditions: Each
$T\in\Gamma^{2}_1(\lambda)$ has a unique major vertex $w_1$ such
that
    \begin{enumerate}[(1)]
   \item  $T - w_1$ has  exactly three components all  from
$\Gamma^{2}_{0}(\lambda)$, where  the neighborhood of $w_1$ in each
component is a pendant vertex of that component;
    \item or   $T - w_1$ has   exactly one component from $\Gamma^{2}_{0}(\lambda)$
    and all other components come from $\Gamma_{0}(\lambda) $, where the neighborhood
    of $w_1$ in each component is a pendant vertex of that component.
    \end{enumerate}

\item  By  $\Gamma^{2}_{2}(\lambda) $, we denote the set of all trees $T$ satisfying
the following conditions:
    Each $T\in\Gamma^{2}_2(\lambda)$ with $\gamma (T ) = 2$ and there is a major vertex $w_2$ of $T$ such that
    \begin{enumerate}[(1)]
    \item  $T - w_2$ has exactly one component from $\Gamma^{2}_{1}(\lambda)$ and all other components come from $\Gamma_{0}(\lambda)$, where  the neighborhood of $w_2$ in each component is a pendant vertex of that component;
    \item or $T - w_2$ has exactly one component, say $T_1$, from $\Gamma_{1}(\lambda)$
    and all other components come from $\Gamma_{0}(\lambda)$, where the unique neighborhood
    of $w_2$ lying in $T_1$ is not a pendant vertex of $T_1$;
    other neighborhood of $w_2$ in each component is a pendant vertex of that component.
    \end{enumerate}

\item  By  $\Gamma^{2}_{3}(\lambda) $, we denote the set of
all trees $T$ satisfying the following three conditions:
    Each $T\in\Gamma^{2}_3(\lambda)$ with $\gamma (T ) = 3$ and there is a major vertex $w_3$
    of $T$ such that
     \begin{enumerate}[(1)]
    \item $T - w_3$ has exactly one component from $\Gamma^{2}_{2}(\lambda)$ and all
    other components come from $\Gamma_{0}(\lambda)$, where the neighborhood of $w_3$
    in each component is a pendant vertex of that component;
   \item  or  $T - w_3$ has exactly one component, say $T_1$, from $\Gamma_{2}(\lambda) $ and all other components come from $\Gamma_{0}(\lambda)$, where  the unique neighborhood of $w_3$ lying in $T_1$ is not a pendant vertex of $T_1$;  other neighborhood of $w_3$ in each component is a pendant vertex of that component;
   \item  or  $T - w_3$ has exactly one component, say $T_1$,
   from $\Gamma_{2}(\lambda) $;  exactly one component, say $T_2$,
   from $\Gamma^{2}_{0}(\lambda) $ and all other components come from
   $\Gamma_{0}(\lambda) $,  where  the neighborhood of $w_3$ in each
   component is a pendant vertex of that component.
   \end{enumerate}

\end{itemize}

Based on the definition for $\Gamma^{2}_{j-1}(\lambda) $, we can
define  $\Gamma^{2}_{j}(\lambda)$, recursively.

 \begin{itemize}

\item By  $\Gamma^{2}_{j}(\lambda) $, we denote the set of all trees $T$ satisfying
the following conditions:
    Each $T\in\Gamma^{2}_j(\lambda)$ with $\gamma (T ) = j$ and there is
    a major vertex $w_j$ of  $T$ such that
    \begin{enumerate}[(1)]
    \item $T - w_j$ has exactly one component from $\Gamma^{2}_{j-1}(\lambda)$ and all other components come from $\Gamma_{0}(\lambda)$, where  the neighborhood of $w_j$ in each component is a pendant vertex of that component;
    \item  or $T - w_j$ has exactly one component, say $T_1$, from $\Gamma_{j-1}(\lambda)$ and all other components come from $\Gamma_{0}(\lambda)$, where  the unique neighborhood of $w_j$ lying in $T_1$ is not a pendant vertex of $T_1$;
    other neighborhood of $w_j$ in each component is a pendant vertex of that component;
    \item  or  $T - w_j$ has exactly one component, say $T_1$, from $\Gamma_{j-1}(\lambda)$;
    exactly one component, say $T_2$, from $\Gamma^{2}_{0}(\lambda)$ and all other
    components come from $\Gamma_{0}(\lambda)$,  where the neighborhood of $w_j$
    in each component is a pendant vertex of that component.
\end{enumerate}
 \end{itemize}

\begin{figure}[h]
\centering
\includegraphics[width=0.8\textwidth]{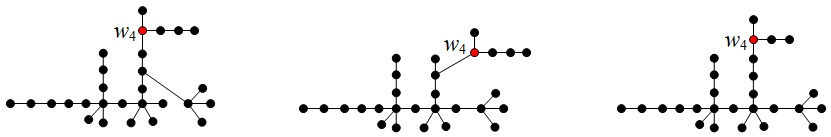}
\caption{ Some trees in $\Gamma^{2}_{4}(0)$}\label{Fig.2}
\end{figure}

\begin{figure}[h]
\centering
\includegraphics[width=0.45\textwidth]{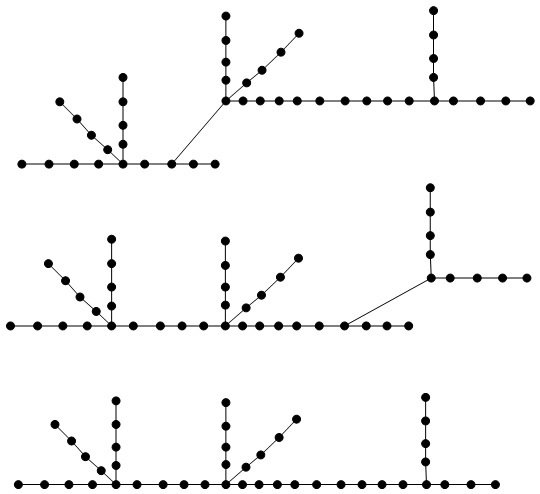}
\caption{ Some trees in   $\Gamma^{2}_{3}(2\cos \frac{\pi}{5})
$}\label{Fig.3}
\end{figure}

Fig.~\ref{Fig.2} provides some trees in $\Gamma^{2}_{4}(0) $, where
$0 = 2\cos \frac{\pi}{2}$. And Fig.~\ref{Fig.3} provides some trees
in $\Gamma^{2}_{3}(2\cos \frac{\pi}{5})$ with three major vertices.

The rest of this paper is organized as follows. In Section 2,  we give some notations and
preliminary lemmas which will be used in our proofs. In Section 3, the proofs for Theorems \ref{1.40} and \ref{1.4} are  presented.

% You may scarsely use \clearpage to advance to a new page if this
% improves the readability of the document structure
\clearpage
\section{Preliminaries}
\label{sec:first}

In this section, we  give some notations and preliminary lemmas which will be used in our proofs.

For a graph $G$, we use $x\sim y$ to denote that $x, y$ are adjacent
vertices (in $G$). In this case, $y$ is called a \emph{neighbor} of
$x$. The \emph{neighborhood} of a vertex $x$ in $G$ is defined as
$N_G (x) = \{y\in V (G)|x \sim y\}$, and the \emph{degree} of $x$ in $G$ is
defined as $d_G( x) = |N_G(x)|$.  In a connected graph $G$ with at
least two vertices, a \emph{pendant vertex} is a vertex of degree 1
and a \emph{quasi-pendant vertex} is the vertex adjacent to a
pendant vertex; a \emph{major vertex} is a vertex with degree at
least 3. By $p(G) $(resp., $\gamma (G$)) we denote the \emph{number
of pendant vertices} in $G$ (resp., major vertices in $G$). If $G$
has only one vertex, we assume $p(G) = 2$ and the single vertex is
also viewed as a pendant vertex of $G$ (without this assumption, the
main result of this article will go wrong). Thus, $p(T ) \geq 2$ for
any tree $T$ and $p(T ) = 2$ if and only if T is a path. By the
\emph{length of a path} we mean the number of edges in the path.
%For $x,y\in V(G)$, the \emph{distance between $x$ and $y$},  denoted by $d(x, y)$,  is the shortest length of paths between $x$ and $y$.
If $S$ is a set of vertices of a graph $G$, we denote by $G - S$ the
(induced) subgraph obtained from $G$ by deleting vertices in $S$
(together with the incident edges). If $S = \{x\}$ or $\{x, y\}$,
then $G - S$ is abbreviated to $G - x$ or $G - x - y$. If $H$ is an
induced subgraph of $G $ with $V(H) \bigcap S = \emptyset$, we use
$H + S$ to denote the subgraph of $G$ induced by the vertex set
$V(H) \bigcup S$.

The following two results are useful in this topic, the first one is
called the Parter-Wiener theorem (see~\cite{Johnson, Parter}
or~\cite{Wiener}) and the second one is due to Johnson.

\begin{lemma}[Parter-Wiener Theorem, \cite{Johnson}]\label{2.1}
 Let $T$ be a tree and suppose that there exists a vertex $v$ of $T$ and a real number $\lambda$ such that $\lambda$ is a common eigenvalue of  $T$ and $T - v$. Then
\begin{enumerate}[(i)]
\item there is a vertex $w$ of  $T$ such that $m(T - w, \lambda) = m(T , \lambda) + 1$;
\item if $m(T , \lambda)\geq2$, then $w$ may be chosen so that $d_T( w) \geq 3$ and so that there are at least three components $T_1$, $T_2$ and $T_3$ of $T - w$ such that $m(T_i, \lambda) \geq 1$, $i = 1, 2, 3$;
\item if $m(T , \lambda) = 1$, then w may be chosen so that $d_T(w )\geq 2$ and so that there are two components $T_1$, $T_2$ of $T -w$ such that $m(T_i, \lambda) = 1$, $i = 1, 2$.
\end{enumerate}
\end{lemma}

\begin{lemma}[\cite{Johnson}]\label{2.2}
Let $T$ be a tree with a vertex $w$. Assume that $\lambda$ is an eigenvalue of $T -w$. Then $m(T - w, \lambda) =m(T , \lambda)+1$ if and only if there is a component $H$ of $T -w$ such that $m(H -v, \lambda) = m(H, \lambda)-1$, where $v$ is the unique neighborhood of $w$ lying in $H$.
\end{lemma}

The following lemma is a proposition of  a tree $T$ in
$\Gamma_k(\lambda)$ which is given in~\cite{Wong}.

\begin{lemma}[\cite{Wong}]\label{2.3}
Let $T \in \Gamma_k(\lambda)$ with $k = \gamma(T )$ and $\lambda=
2\cos \left(\frac{i\pi}{m+1}\right)$, where $i\in [1, m]$ and $m +
1$ are co-prime integers. Then the following assertions hold.
\begin{enumerate}[(i)]
\item $\lambda$ is an eigenvalue of $T$;

\item  If $v$ is a pendant vertex of $T$, then $m(T - v, \lambda) = m(T , \lambda) - 1$.
 \end{enumerate}
\end{lemma}

\begin{lemma}[\cite{Johnson}]\label{2.4}
Let $P_n$ be  a path on $n$ vertices.  Then for any eigenvalue
$\lambda$ of $T$, we have $m(P_n,\lambda)=1$.
\end{lemma}

\section{Proofs of Theorems \ref{1.40} and \ref{1.4}}
\label{sec:hints}
{\bf Proof of Theorem~\ref{1.40}:}\ Let $T$ be a tree with $p(T ) =3
$. Since $p(T ) =3$, $T$ has exactly one major vertex, say $w_1$,
such that  $T -w_1 = T_1 \bigcup T_2 \bigcup  T_3 $, where $T_1,T_2,
T_3$  all are the path-components of $T -w_1$. Let $u_i$ be the
neighborhood of $w_1$ in $T_i$ ($i=1,2,3$).

``If" part: Suppose that $T \in\Gamma^{2} _1(\lambda)$, by the definition
of  $\Gamma^{2}_{1}(\lambda) $, we consider the following two conditions.

\begin{enumerate}[(i)]

\item If  $T_i\in \Gamma^{2}_{0}(\lambda) $ for  $i=1,2,3$. Then $m(T
-w_1 ,\lambda)=0.$ By Interlacing Theorem (see~\cite{Brouwer}), $m(T
,\lambda)\leq m(T -w_1 ,\lambda)+1=1.$ Since $\lambda$ is  an
eigenvalue of $T$, $m(T ,\lambda)\geq 1$. Hence, $ m(T,\lambda)=
1=p(T)-2$, as desired.

\item If there is exactly one component of $T-w_1$, say $T_1$, which
is from $\Gamma^{2}_{0}(\lambda) $ and $T_2,T_3$ come from
$\Gamma_{0}(\lambda) $. Then, $m(T_1,\lambda)=0$. By
Theorem~\ref{1.3}, $m(T_j,\lambda)=p(T_j)-1=1$ $(j=2,3)$ and
$m(T_j-u_j,\lambda)=m(T_j,\lambda)-1=0$. Applying Lemma~\ref{2.2},
we have
$$m(T , \lambda)=-1+m(T-w_1 , \lambda) = -1+0+1+1=1.$$
Hence, $m(T,\lambda)= 1=p(T)-2$, as desired.

\end{enumerate}

``Only if" part: Suppose that $\lambda$ as an eigenvalue of  $T$
with multiplicity $p(T ) -2=1$. We consider the following two cases.

\textbf{Case 1.} $m(T-w_1,\lambda)\geq1$.

Hence, $\lambda$ is also an eigenvalue of $T -w_1$. Among these
path-components, at least one component, say $T_3$, has $\lambda$ as
an eigenvalue. Then Lemma~\ref{2.4} implies that $m(T_3,\lambda)=1$.
Hence, $T_3\in \Gamma_0(\lambda)$. By Lemma~\ref{2.3},
$m(T_3-u_1,\lambda)=0=m(T_3,\lambda)-1$. Applying Lemma~\ref{2.2},
we then have
$$m(T-w_1 , \lambda) =  m(T , \lambda) +1=2.$$
Hence, there must be a path-component,  say $T_2$, has $\lambda$ as
an eigenvalue. But the other path-component $T_1$, does   not  have
$\lambda$ as an eigenvalue.  Then  $T_2\in \Gamma_0(\lambda)$ and
$T_1\in \Gamma^{2}_0(\lambda)$. It means that there is a major
vertex $w_1$ of  $T$ such that $T - w_1$ has  exactly one component
from $\Gamma^{2}_{0}(\lambda)$ and all other components come from
$\Gamma_{0}(\lambda)$, where  the neighborhood of $w_1$ in each
component is a pendant vertex of that component. Thus,
$T\in\Gamma^{2}_{1}(\lambda)$ .

\textbf{Case 2.}  $m(T-w_1,\lambda)=0$.

In this case, $m(T_i,\lambda)=0$, that is, $T_i\in \Gamma^{2}_0(\lambda)$ for $i=1,2,3.$
 It means that there is a major vertex $w_1$ of  $T$ such that $T - w_1$ has  exactly three components all  from $\Gamma^{2}_{0}(\lambda) $, where  the neighborhood of $w_1$ in each component is a pendant vertex of that component. Thus, $T\in\Gamma^{2}_{1}(\lambda) $ .

The proof is completed.\q

{\bf Proof of Theorem~\ref{1.4}:}\ ``If" part: Clearly, $\gamma (T )
\geq 1$ since $p (T ) \geq 4$. We proceed by induction on $\gamma (T
)$ to prove that $m(T , \lambda) = p(T ) - 2$. If $\gamma (T ) = 1$,
then $T \in \Gamma^{2}_{1}(\lambda)$.   Thus, $T$ has exactly one
major vertex, say $w_1$, such that  $T -w_1= T_1 \bigcup T_2 \bigcup
\cdots  \bigcup T_s$, where  $s\geq4$ and  $T_1,\ldots, T_s$ all are
the path-components of $T -w_1$. Let $u_i$ be the neighborhood of
$w_1$ in $T_i$ ($i=1,\ldots,s$).

 By definition of  $\Gamma^{2}_{1}(\lambda) $ and $p (T ) \geq 4$, there is
 exactly one component, say
$T_1$, which is from $\Gamma^{2}_{0}(\lambda) $ and $T_2,\ldots,T_s$
come from $\Gamma_{0}(\lambda) $, where   each $u_i$ is  a pendant
vertex of $T_i$ $(i=1,\ldots,s)$. Then, $m(T_1,\lambda)=0$. By
Theorem \ref{1.3}, $m(T_j,\lambda)=p(T_j)-1=1$ $(j=2,\ldots,s)$ and
$m(T_j-u_j,\lambda)=m(T_j,\lambda)-1=0$. Since $\lambda$ is also an
eigenvalue of $T -w_1$, by Lemma \ref{2.2},
$$ m(T , \lambda)=-1+m(T-w_1 , \lambda) = -1+0+s-1=s-2.$$
Hence, $ m(T,\lambda)= s-2=p(T)-2$.

Suppose that $m(T , \lambda) = p(T ) -2$ for all trees $T$ in
$\Gamma^{2}_{j}(\lambda) $ with $j \leq k -1$ and $T$ lies in
$\Gamma^{2}_{k}(\lambda)$, where  $\gamma(T )=k\geq 2$.  By the
definition of  $\Gamma^{2}_{k}(\lambda) $,  there is a major vertex
$w_k$ such that $T -w_k= T_1 \bigcup T_2 \bigcup \cdots  \bigcup
T_s$,  where $s\geq3$. Let $u_i$ be the  neighborhood of $w_k$ in
$T_i$  ($i=1,\ldots,s$).

 By the construction of $\Gamma^{2}_{k}(\lambda)$, we now consider the following three conditions.

\begin{enumerate}[(i)]
\item If there is  exactly one component, say  $T_1$, which is from
$\Gamma^{2}_{k-1}(\lambda) $ and  $T_2,\ldots,T_s$ come from
$\Gamma_{0}(\lambda) $, where  $u_i$ is  a pendant vertex of $T_i$
$(i=1,\ldots,s)$. Since $\gamma (T_1) < k$, the induction hypothesis
implies that $m(T_1, \lambda) = p(T_1) -2$. By Lemma \ref{1.3},
$m(T_j,\lambda)=p(T_j)-1=1$ $(j=2,\ldots,s)$ since  $T_j\in
\Gamma_{0}(\lambda)$. Then we have
\begin{align*}
m(T-w_k,\lambda)&= \sum_{i=1}^{s}m(T_i,\lambda)\\
&=p(T_1)-2+(s-1)\\
&=p(T)-1.
\end{align*}
The last equality holds from    $u_i$ $(i=1,\ldots,s)$ is a pendant
vertex of $T_i$ but not that of $T$, we have $p(T ) = (p(T_1)-1) +
s-1$.

\item If there is  exactly one component, say $T_1$, which is from
$\Gamma_{k-1}(\lambda) $ and $T_2,\ldots,T_s$ come from
$\Gamma_{0}(\lambda) $, where  $u_1$ is not a pendant vertex of
$T_1$;   $u_j$   is a pendant vertex of $T_j$ $(j=2,\ldots,s)$. By
Theorem \ref{1.3},  $m(T_i,\lambda)=p(T_i)-1$ $(i=1,\ldots,s)$. Then
we have
\begin{align*}
m(T-w_k,\lambda)&= \sum_{i=1}^{s}m(T_i,\lambda)\\
&= \sum_{i=1}^{s}(p(T_i)-1)\\
&=\sum_{i=1}^{s}p(T_i)-s\\
&=p(T)-1.
\end{align*}
The last equality holds from $u_1$ is not a pendant vertex of $T_1$;
$u_j$ ($j=2,\ldots,s$) is a pendant vertex of $T_j$ but not that of
$T$, we have $p(T ) = \sum_{i=1}^{s}p(T_i)- s+1$.

\item If there is  one component, say  $T_1$, which is
 from $\Gamma_{k-1}(\lambda)$;  one component, say $T_2$, which  is
 from $\Gamma^{2}_{0}(\lambda) $ and $T_3,\ldots,T_s$ come
 from $\Gamma_{0}(\lambda)$, where  $u_i$ is  a pendant
 vertex of $T_i$ $(i=1,\ldots,s)$. Then, we have $m(T_2,\lambda)=p(T_2)-2=0$.
 By Lemma \ref{1.3}, $m(T_1,\lambda)=p(T_1)-1$ and  $m(T_l,\lambda)=p(T_l)-1=1$ $(l=3,\ldots,s)$.
Then we have
\begin{align*}
m(T-w_k,\lambda)&= \sum_{i=1}^{s}m(T_i,\lambda)\\
&=(p(T_1)-1) +(p(T_2)-2) +\sum_{l=3}^{s}(p(T_l)-1)\\
&=\sum_{i=1}^{s}p(T_i)-s-1\\
&=p(T)-1.
\end{align*}
 The last equality holds from    $u_i$ $(i=1,\ldots,s)$ is a pendant vertex of $T_i$ but not that of $T$, we have $p(T ) = \sum_{i=1}^{s}p(T_i)- s$.
\end{enumerate}

From the above consideration, we have $m(T-w_k,\lambda)=p(T)-1$.
 By the Interlacing Theorem, we have $m(T,\lambda)\geq m(T-w,\lambda)-1=(p(T)-1)-1=p(T)-2.$  By Theorem \ref{1.1} and Theorem \ref{1.3},
 $ m(T,\lambda)\leq p(T)-2$. Hence, $ m(T,\lambda)= p(T)-2$.

``Only if" part: Suppose that $\lambda$ is an eigenvalue of $T$ with
multiplicity $p(T ) -2$. We first prove the following  three Claims.

{\bf Claim 1.} For any $v\in V(T)$, $ m(T-v,\lambda)\geq 1$.

 Assume to contrary that there is  a vertex  $v\in V(T)$ such that $ m(T-v,\lambda)=0$.  By the Interlacing Theorem,
we have $ m(T,\lambda)\leq m(T-v,\lambda)+1= 1$, which contradicts with $ m(T,\lambda)=p(T)-2\geq 2$ since $p(T)\geq 4$. \q

There must be a major vertex $w$, $d_T( w) = s \geq 3$, such that $T
- w$ has at least $s -1$ path-components. Suppose $T -w = T_1
\bigcup T_2 \bigcup \cdots \bigcup T_s $, where $T_2,\ldots, T_s$
are the path-components of $T -w$. Let $u_i$ be the unique
neighborhood of $w$ in $T_i$ $(i=1,\ldots,s)$.

 {\bf Claim 2.}  $m(T_i, \lambda) \geq p(T_i) - 2$ ($i=1,\ldots,s$).

We assume to contrary that there is a component of $T-w$, say $T_1$, such that $m(T_1, \lambda) \leq p(T_1) - 3$. By Claim 1, $\lambda$ is also an eigenvalue of $T -v$ for any vertex $v$ of $T$. Applying Lemma \ref{2.1}, among these path-components, at least one component, say $T_s$, has $\lambda$ as an eigenvalue. Then $\lambda$ is not an eigenvalue of  $T_s -u_s $.  By Lemma \ref{2.2}, we have
$$p(T) - 2=m(T , \lambda) = -1 + m(T - w, \lambda) .$$
Then we have $m(T - w, \lambda)= p(T) - 1.$
On the other hand,
\begin{align*}
m(T - w, \lambda)&= \sum_{i=1}^{s}m(T_i,\lambda)\\
&=m(T_1,\lambda)+\sum_{j=2}^{s}m(T_j,\lambda)\\
&\leq(p(T_1)-3) +\sum_{j=2}^s (p(T_j)-1)\\
&=\sum_{i=1}^{s}p(T_i)-s-2\\
&\leq p(T)-2,
\end{align*}
which contradict with $m(T - w, \lambda)= p(T) - 1.$\q

  {\bf Claim 3.}  If there is one component of $T-w$, say $T_k$,
  which satisfies $m(T_k, \lambda)= p(T_k) - 2$, then other $T_j$ ($j\neq k$) satisfies $m(T_j, \lambda)= p(T_j) - 1$.

 By Claim 2, $m(T_i, \lambda) \geq p(T_i) - 2$ ($i=1,\ldots,s$).
 Assume to contrary that there are at least two  components of $T-w$, say $T_1$ and $T_2$, such that $m(T_1, \lambda) =p(T_1) - 2$ and $m(T_2, \lambda) =p(T_2) - 2$. Similar with the proof of Claim 2, we have
\begin{align*}
m(T - w, \lambda)&= \sum_{i=1}^{s}m(T_i,\lambda)\\
&=m(T_1,\lambda)+m(T_2,\lambda)+\sum_{j=3}^{s}m(T_j,\lambda)\\
&\leq(p(T_1)-2) +(p(T_2)-2)+\sum_{j=3}^s (p(T_j)-1)\\
&=\sum_{i=1}^{s}p(T_i)-s-2\\
&\leq p(T)-2,
\end{align*}
which contradict with $m(T - w, \lambda)= p(T) - 1.$\q

By Claim 1, $\lambda$ is also an eigenvalue of $T -v$ for any vertex $v$ of $T$. Among these path-components, at least one component, say $T_s$, has $\lambda$ as an eigenvalue. Then $\lambda$ is not an eigenvalue of  $T_s -u_s $.  Applying Lemma \ref{2.2},
$$m(T , \lambda) = -1 + m(T - w, \lambda) = -1 +\sum_{i=1}^s m(T_i, \lambda).$$
Since $T_2,\ldots, T_s$ are the path-components of $T -w$, $u_j$ is
a pendant vertex of $T_j$ $(j=2,\ldots,s)$. According to the degree
of $u_1$ in $T_1$, we now consider the following two cases.

\textbf{Case 1.}   $u_1$ is not a pendant vertex of $T_1$.

Since $m(T_i,\lambda) \leq p(T_i) -1$ for each $T_i $
$(i=1,\ldots,s)$, it follows from the above equality, we have
$$m(T , \lambda) \leq-1 +\sum_{i=1}^s (p(T_i)-1)=-1+\sum_{i=1}^s p(T_i)-s.$$
Clearly, $p(T ) = \sum_{i=1}^s p(T_i)-s+1$ since $u_1$ is not a
pendant vertex of $T_1$  and  $u_j$ ($j=2,\ldots,s$) is a pendant
vertex of $T_j$ but not that of $T$. Then we have
$$p(T ) -2 = m(T , \lambda)\leq -1  +\sum_{i=1}^s p(T_i)-s= p(T ) - 2.$$
Thus, $m(T_i, \lambda) = p(T_i) - 1$ for each $i=1,\ldots,s$. Then
each $T_j$ ($j=2,\ldots,s$) lies in $\Gamma_{0}(\lambda) $ and $T_1$
lies in $\Gamma_{k-1}(\lambda) $ since $\gamma (T_1) = k-1$  (by
Theorem \ref{1.3}). Hence, we have $T \in\Gamma^{2}_{k}(\lambda)$.

\textbf{Case 2.}   $u_1$ is a pendant vertex of $T_1$.

We  claim that there must be a component, say $T_t$, satisfies
$m(T_t, \lambda) = p(T_t) - 2$. Otherwise, by Claim 2, $m(T_i,
\lambda) = p(T_i) - 1$ for $i=1,\ldots,s$. Then $T - w$ has exactly
one component $T_1$ from $\Gamma_{k-1}(\lambda) $ and all other
components $T_j$ $(j=2,\ldots,s)$ come from $\Gamma_{0}(\lambda) $,
where  $u_i$ is a pendant vertex of $T_i$  $(i=1,\ldots,s)$. It
follows from   Theorem \ref{1.3}, we have   $m(T, \lambda) = p(T) -
1$ which contradict with the assumption.

If $t=1$, then $m(T_1, \lambda) = p(T_1) - 2$. By Claim 3, $m(T_j,
\lambda) = p(T_j) - 1$ for $j=2,\ldots,s$.
Clearly, $\gamma (T ) \geq 1$ since $p (T ) \geq 4$. We proceed by
induction on $\gamma (T )$ to prove that
$T\in\Gamma^{2}_{k}(\lambda) $, where $k=\gamma (T )$.  If $\gamma
(T ) = 1$, then  $T$ has exactly one major vertex, say $w_1$, such
that  $T -w_1= T_1 \bigcup T_2 \bigcup \cdots  \bigcup T_s$, where
$s\geq4$ and  $T_1,\ldots, T_s$  all are the path-components of $T
-w_1$. By $m(T_1, \lambda) = p(T_1) - 2$ and $m(T_j, \lambda) =
p(T_j) - 1$ ($j=2,\ldots,s$), we have $T_1$ is from
$\Gamma^{2}_{0}(\lambda) $ and $T_2,\ldots,T_s$ come from
$\Gamma_{0}(\lambda) $ (by Theorem \ref{1.3}), where  each $u_i$ is
a pendant vertex of $T_i$ $(i=1,\ldots,s)$. Then by definition of
$\Gamma^{2}_{1}(\lambda)$, $T\in\Gamma^{2}_{1}(\lambda)$.
Suppose that $T\in\Gamma^{2}_{j}(\lambda)$ for all trees with $m(T ,
\lambda) = p(T ) -2$ where  $j =\gamma (T )\leq k -1$ and $\gamma (T
)=k\geq 2$. Since $\gamma (T_1) =k-1<k$, the induction hypothesis
implies that $T_1\in\Gamma^{2}_{k-1}(\lambda)$.  By Theorem
\ref{1.3}, $T_j\in \Gamma_{0}(\lambda)$ $(j=2,\ldots,s)$. Then by
definition of $\Gamma^{2}_{k}(\lambda)$,
$T\in\Gamma^{2}_{k}(\lambda)$.

If $t\in\{2,\ldots,s\}$,  without loss of generality,  we can assume
that $t=2$. Then $m(T_2, \lambda) = p(T_2) - 2$. By Claim 3, $m(T_j,
\lambda) = p(T_j) - 1$ for $j\neq2$ and $j\in\{1,\ldots,s\}$. By
Theorem \ref{1.3},   $T_1\in \Gamma_{k-1}(\lambda)$ and $T_j\in
\Gamma_{0}(\lambda)$ $(j=3,\ldots,s)$. $T_2\in
\Gamma^{2}_{0}(\lambda)$ since  $m(T_2, \lambda) = p(T_2) - 2$. Then
by the definition of $\Gamma^{2}_{k}(\lambda)$, we have
$T\in\Gamma^{2}_{k}(\lambda)$. This completes the proof of
Theorem~\ref{1.4}. \q

\section{Concluding remarks}

In this paper, we study the multiplicity of an eigenvalue $\lambda$ of a  tree and characterize all trees with eigenvalue multiplicity two less than their number of pendant vertices, that is
 $m(T,\lambda)=p(T)-2$. But it seems somewhat difficult to give a characterization
 on trees $T$ with the eigenvalue $\lambda$ such
that $m(T, \lambda) = p(T) - k$ for any $k\in [1, p(T)-1]$. We leave
it for further study.

\acknowledgements
\label{sec:ack}
The authors would like to thank the anonymous referees for their constructive corrections and valuable
comments on this paper, which have considerably improved the presentation of this paper.

\nocite{*}
\bibliographystyle{abbrvnat}
% use the following instead if you encounter problems
%\bibliographystyle{alpha}
\bibliography{sample-dmtcs}
\label{sec:biblio}

\end{document}